\documentclass{article}
\usepackage{latexsym}
\usepackage{amssymb}
\usepackage{amsmath}
\usepackage{color}
\usepackage{graphicx}
\allowdisplaybreaks

\newtheorem{theorem}{\color{black} Theorem}[section]
\newtheorem{lemma}{\color{black} Lemma}[section]

\begin{document}

\title{\bf Persistence of Hyperbolic Tori in Generalized Hamiltonian
Systems\thanks{This paper has been published in: {\it Northeast.
Math. J.} {\bf 21} (4) (2005), 447-464.}}

\author{
Zhenxin Liu\thanks{{Corresponding author:} zxliu@email.jlu.edu.cn
(Zhenxin Liu).}, Dalai Yihe and Qingdao Huang
\\{\small College of Mathematics, Jilin University,
Changchun 130012, P. R. China} }
\date{}
\maketitle {\noindent{\bf Abstract}}

In this paper we prove the persistence of hyperbolic invariant tori
in generalized Hamiltonian systems, which may admit a distinct
number of action and angle variables. The systems under
consideration can be odd dimensional in tangent direction. Our
results generalize the well-known results of Graff and Zehnder in
standard Hamiltonians. In our case the unperturbed Hamiltonian
systems may be degenerate. We also consider the persistence
problem of hyperbolic tori on sub-manifolds.

\vskip 5mm
\section {Introduction and Main Result}
According to the celebrated KAM (Kolmogorov-Arnold-Moser) theory, we
know that most of invariant tori of integrable Hamiltonians persist
under a small perturbation. In their case, the Hamiltonian is
standard and the dimension of invariant tori equals the degree of
freedom, i.e., the ``highest dimensional tori''. However, some
Hamiltonian systems' highest dimensional tori cannot survive the
perturbations, but some tori which have lower dimension can be
persisted under small perturbations, which are called lower
dimensional invariant tori. In 1965, Melnikov$^{[1]}$ formulated a
KAM type persistence result for elliptic lower dimensional tori of
integrable Hamiltonian systems under so-called Melnikov's
non-resonance condition. But the complete proof of his result was
carried out more than twenty years later by Eliasson, Kuksin, and
P\"{o}schel (see [2]--[4]). For the persistence of hyperbolic lower
dimensional tori, in 1974, Graff$^{[5]}$ considered the following
Hamiltonian system:
\begin{equation}\label{1.11}
H=e+\langle \omega_0,y\rangle+\frac12\langle y,Ay\rangle+
\frac12\langle z,Mz\rangle+P(x,y,z),
\end{equation}
where $(x, y, z)\in T^{n}\times R^{n}\times R^{2m}$,
$M=\left(\begin{array}{cc}O &B\\B^{\top}&O\end{array}\right)$,
$\omega_0\in R^n$ is a fixed Diophantine toral frequency, and $P$ is
a small perturbation. The persistence of the unperturbed Diophantine
hyperbolic torus $T^n\times\{0\}\times\{0\}$ was shown as well as
the preservation of the toral frequency $\omega_0$. Zehnder in [6],
using generalized implicit function theorem, proved the same result.
More recently, Li and Yi$^{[7]}$ generalized the results of Graff
and Zehnder on the persistence of hyperbolic invariant tori in
Hamiltonian systems by allowing the degeneracy of the unperturbed
Hamiltonians and they obtain the preservation of part or full
components of frequencies. They adopted the Fourier series expansion
for normal form $N$, which is a new technique.

Due to important technical reasons, the development of KAM theory
for odd dimensional systems has been considered as a challenging
problem. Li and Yi$^{[8]}$  solve the delicate problem by
considering generalized Hamiltonian systems which preserve  a
prescribed Poisson structure instead of volume. In their case, the
Hamiltonians considered may admit distinct number of action and
angle variables and more important, which can be odd dimensional.
Motivated by their work, in this paper, we show that the
Graff-Zehnder result also holds in generalized Hamiltonian systems.

We consider the following parameter-dependent Hamiltonian system:
  \begin{equation}\label{1.1}
H=e(\lambda)+\langle \Omega(\lambda),y\rangle+\frac12\left\langle
{y\choose z},{\mathcal M} (x,\lambda){y\choose
z}\right\rangle+h(x,y,z,\lambda)+ P(x,y,z,\lambda),
\end{equation}
where $(x, y, z)\in T^{n}\times R^{l}\times R^{2m}$, $\lambda$ is a
parameter in a bounded, closed, connected region ${\Lambda}\subset
R^k$, ${\mathcal M}$ is  symmetric, real analytic in $x\in {\mathcal
D}(r)=\{x\in {C^n}/{Z^n}:\ |{\rm Im}x|<r\}$,
$h(x,y,z,\lambda)=O(|(y,z)|^3)$ is real analytic, and, the
perturbation $P$ is real analytic in a complex neighborhood $D(r,
s)=\{(x,y,z) : \ |{\rm Im} x|<r, |y|<s, |z|<s\}$ of
$T^n\times\{0\}\times\{0\}$. In the above, all $\lambda$ dependence
are of class $C^{l_0}$ for some $l_0\ge n$.

Write ${\mathcal M}$ in (\ref{1.1}) into blocks:
\begin{equation}\label{block}
{\mathcal M}=\left(\begin{array}{cc}A
&B\\B^\top&M\end{array}\right),
\end{equation}
where $A=A(x,\lambda)$, $B=B(x,\lambda)$, $M=M(x,\lambda)$ are
$l\times l$, $l\times 2m$, $2m\times 2m$ minors of $\mathcal
M={\mathcal M}(x,\lambda)$ respectively.

A so-called generalized Hamiltonian system is defined on a Poisson
manifold which can be odd dimensional and structurally degenerate.
Consider the manifold $G\times T^n\times R^{2m}$, where $G\subset
R^l$ is a bounded, connected and closed region, $T^n$ is the
standard $n$-torus and $l,n,m$ are positive integers. Let $I$ be the
structure matrix in tangent direction, and $J$ be the $2m\times 2m$
standard symplectic matrix in norm direction. As in [8], assume
$I=I(\lambda)$ be real analytic. Then the structure matrix $\tilde
I$ on $G\times T^n\times R^{2m}$ has the following form:
$$\tilde I(\lambda)=\left(\begin{array}{cc}I(\lambda) &O\\
O&J\end{array}\right),~~~~
I(\lambda)=\left(\begin{array}{cc}O &E(\lambda)\\
-E^\top(\lambda)&C(\lambda)\end{array}\right),$$ where $O$ denotes
zero matrix, $E=E_{l,n},C=C_{n,n}$ with $C^\top=-C$. Let $\nabla$
denote the standard Euclidean gradient on $R^l\times T^n\times
R^{2m}$. Then $\tilde I$ defines a Poisson structure or a 2-form
$\omega^2$ in the following way:
$$\{f_1,f_2 \}=df_2(\tilde I df_1)=\langle \nabla f_1,\tilde I\nabla f_2 \rangle
=\omega^2(\tilde I df_1,\tilde I df_2),$$ for  all smooth functions
$f_1,f_2$ defined on $G\times T^n\times R^{2m}$, where
$\{\cdot,\cdot\}$ denotes the Poisson bracket.

Then the equation of motion associated with (\ref{1.1}) reads
$$
\left(\begin{array}{c}\dot y\\ \dot x\\ \dot z
\end{array}\right)=\tilde I \nabla H.
$$
Thus, the unperturbed system associated with (\ref{1.1}) admits a
smooth family of invariant $n$-tori  $T_{\lambda}=T^n\times
\{0\}\times\{0\}$ with toral frequencies
$\omega(\lambda)=-E^\top(\lambda)\Omega(\lambda)$ parameterized by
$\lambda\in \Lambda$. As in [7], we first assume that $J[M]$ is
hyperbolic on $\Lambda$, i.e., if $\lambda_i(\lambda)$,
$i=1,2,\cdots,2m$, are
  eigenvalues
of $J[M](\lambda)$, then

{\bf H)} there exists a constant $\sigma_0>0$ such that
$$
|{\rm Re}\lambda_i(\lambda)|\ge\sigma_0,
$$
for all $\lambda\in \Lambda$ and $i=1,\cdots,2m$.

Next,  we assume the R\"ussmann condition

{\bf R)}
$$
\max_{\lambda\in \Lambda}{\rm
rank}\{\partial^{\alpha}\omega(\lambda): ~~\forall |\alpha|\le
n-1\}=n.
$$

To make a difference between $\Omega(\lambda)$ and toral frequency
$\omega(\lambda)= -E^\top(\lambda)\Omega(\lambda)$, we call
$\Omega(\lambda)$ as pseudo-frequency. For general structure matrix
$I$, we cannot obtain the persistence of part frequency components
by the associate persistence of part pseudo-frequency components.
But for some special structure matrix $I$, we can even obtain the
unchanged toral frequency in spite that only part pseudo-frequency
components are preserved (see Example 5.2), which of course depends
closely on the specific form of $E(\lambda)$ and $\Omega(\lambda)$.
So it is necessary to study the preservation of part or full toral
pseudo-frequency components in connection with the  degree of
non-degeneracy of the matrix $[A]$. As in [7], we assume that

{\bf ND)} there is a $1\le n_0\le l$ such that both the $n_0\times
n_0$ ordered principal minor $U$ of $[A]$ and  $Y\equiv [M]-[B]^\top
{\rm diag}(U^{-1},O)[B]$ are non-singular on $\Lambda$, where $O$
denotes the zero matrix.

It is clear that ND) holds automatically if $[A]$ is non-singular on
$\Lambda$ and $|[B]|_{\Lambda}$ is sufficiently small (in
particular, when $[B]\equiv 0$).

Define
\begin{equation}\label{eta}
\eta=\frac 2{\sqrt{{\rho}_0^2+4\alpha{\rho}_0}+{\rho}_0},
\end{equation}
where
\begin{eqnarray}\label{alpha}
\alpha&=&(1+ 2m)(|Y^{-1}|+|U^{-1}|+(|Y^{-1}||U^{-1}|)
(2|[B]|+|[B]|^2|U^{-1}|))_{\Lambda},\label{alpha}\\
{\rho}_0&=&\frac {4m}{\sigma_0}\left(1+\frac {2m}{\sigma_0}
|[M]|_{\Lambda}\right)^{2m-1}.\label{rho0}
\end{eqnarray}

The main result of this paper is the following.

\begin{theorem} Consider {\rm (\ref{1.1})} and assume the conditions
 {\rm H)}, {\rm R)}, {\rm ND)} and
 \begin{equation}\label{mb}
|M-[M]|_{{\mathcal D}(r)\times \Lambda},~~~|B-[B]|_{{\mathcal
D}(r)\times \Lambda}<\eta.
\end{equation}
Then  there is an  $\varepsilon=\varepsilon(r,s,l_0,\sigma_0,U)>0$
sufficiently small such that if
\begin{equation}\label{1.3}
|\partial^{l}_{\lambda} P|_{D(r,s)\times\Lambda} <
\gamma^{n+1}s^2\varepsilon,~~~~|l|\le l_0,
\end{equation}
then

1) there is a  $0<r_0=r_0(r,\sigma_0,U)\le r$ and a Cantor-like set
$\Lambda_{\gamma}\subset \Lambda$, with $|\Lambda\setminus
\Lambda_{\gamma}| =O(\gamma^{\frac1{n_*-1}})$, where
$n_*=\max\{2,n\}$, for which there is a $C^{l_0-1}$ Whitney smooth
family of real analytic, symplectic transformations
$$
\Psi_{\lambda}: D\left(\frac{r_0}2,\frac{s}{2}\right)\rightarrow
D(r_0,s), \quad \lambda\in \Lambda_{\gamma},
$$
which are  $C^{l_0}$ uniformly close to the identity such that
$$
H\circ \Psi_{\lambda}=e_*+\langle\Omega_*(\lambda),y\rangle+{1\over
2}\left\langle {y\choose z},{\mathcal M}_*(x,\lambda){y\choose
z}\right\rangle+h(x,y,z,\lambda) +P_*(x,y,z,\lambda),
$$
where
\begin{eqnarray*}
&&|\partial^{l}_{\lambda}e_*-\partial^{l}_{\lambda}e|_{\Lambda_{\gamma}}
=O(\gamma^{n+1}s\varepsilon\zeta),\\
&&|\partial^{l}_{\lambda}\Omega_*-\partial^{l}_{\lambda}\Omega|_{\Lambda_{\gamma}}=
O(\gamma^{n+1}s\varepsilon\zeta),\\
&&|\partial^{l}_{\lambda}{\mathcal M}_*-
\partial^{l}_{\lambda}{\mathcal M}|_{{\mathcal D}(r_0)\times
\Lambda_{\gamma}}=O(\gamma^{n+1}\varepsilon\zeta).
\end{eqnarray*}
Thus,  all unperturbed tori $T_{\lambda}$ with
$\lambda\in\Lambda_{\gamma}$ will persist and give rise to a
$C^{l_0-1}$ Whitney smooth family of slightly deformed, analytic,
quasi-periodic, invariant $n$-tori of the perturbed system;

2)
$$
(\Omega_*(\lambda))_i=(\Omega_0(\lambda))_i,~~~\lambda\in\Lambda_{\gamma},~
i=1,2,\cdots, n_0,
$$
i.e., the first $n_0$ components of the perturbed toral
pseudo-frequency remain unchanged. In particular, if $n_0=l$, i.e.,
$U=[A]$ is non-singular on $\Lambda$, then every  Diophantine tori
$T_{\lambda}$ with Diophantine type  $(\gamma,\tau)$ for a fixed
$\tau>n-1$ will persist  with unchanged toral frequencies.
\end{theorem}

\section{KAM Step}
\setcounter{equation}{0}

In this section, we describe the linear iterative scheme with
respect to (\ref{1.1}) for one KAM step, say, form a $\nu$th step to
the $(\nu+1)$th step. Below, let $\tau>\max\{n(n-1)-1,l(l-1)-1,0\}$
be fixed.

Consider (\ref{1.1}) and define $e_0=e$, $\Omega_0=\Omega$,
$\mathcal M^0= \mathcal M$, $A^0=A$, $B^0=B$, $M^0=M$, $h_0=h$,
$P_0=P$, ${\Lambda}_0=\Lambda$, $\gamma_0=\gamma$, $r_*=r$,
$s_0=\left(\frac{\gamma_0}2\right)^{n+1}\varepsilon_0^{ \frac59}$.
We rewrite $[A^0](=[A])$ into blocks:
\[
[A^0]=\left(\begin{array}{ll}U^0&D^0\\(D^0)^\top&V^0\end{array}\right),
\]
where $U^0=U$. Without loss of generality, assume that
$0<s_0,r_0,\varepsilon_0 \le 1$. By (\ref{1.3}), we have
\begin{equation}\label{P0}
|\partial_\lambda^{l} P_0|_{D(r_0,s_0)}\le \gamma_0^{n+1}
s_0^2\varepsilon_0,~~~~|l|\le n.
\end{equation}

In what follows, quantities (domains, normal form, perturbation,
etc.) without subscripts denotes the Hamiltonian in $\nu$-th step,
while those with subscript ``+" denotes the Hamiltonian of
$(\nu+1)$-th step. And we shall use ``$<\cdot$" to denote ``$<c$"
with a constant $c$ which is independent of the iteration step. For
simplicity, we set $l_0=n$.

Suppose that at the $\nu$-th step, we have arrived at the following
Hamiltonian:
\begin{eqnarray}
&&H=N+P, \label{2.1}\\
& &N=e+\langle\Omega(\lambda),y\rangle+\frac12\left\langle {y\choose
z},{\mathcal M}(x,\lambda) {y\choose
z}\right\rangle+h_0(x,y,z,\lambda),\nonumber
\end{eqnarray}
where $(x,y,z)\in D=D(r,s)$, $\lambda\in \Lambda$, $e(\lambda),
\Omega(\lambda)$ are  smooth on $\Lambda$, ${\mathcal
M}(x,\lambda)=\left(\begin{array}{ll}A&
B\\B^\top&M\end{array}\right)$ is real symmetric over ${\mathcal
D}\times {\Lambda}=\{x :|{\rm Im}x|<r\}\times {\Lambda}$ which is
smooth in $\lambda\in \Lambda$ and real analytic in $x\in {\mathcal
D}={\mathcal D}(r)$, $P$ is real analytic in $(x,y,z)\in D$, smooth
in $\lambda\in \Lambda$, and moreover,
$$
|\partial_\lambda^{l} P|_{D\times\Lambda}\le \gamma^{n+1}
s^2\varepsilon,~~~|l|\le n.
$$

We shall construct a symplectic transformation $\Phi=\Phi_{\nu+1}$
which transforms the Hamiltonian (\ref{2.1}), in smaller phase and
frequency domains, to the desired Hamiltonian in the next KAM cycle
(the $(\nu+1)$th KAM step).

Define
\begin{eqnarray}
\varepsilon_+&=&\varepsilon^{\frac{10}9},\nonumber\\
\gamma_+&=&\frac{\gamma_0}4+\frac{\gamma}2,\nonumber\\
r_+&=&\frac{r_0}4+\frac{r}2,\nonumber\\
s_+&=&\frac18\alpha s,\ \alpha=\varepsilon^{\frac13},\nonumber\\
K_+&=&([\log\frac{1}{s}]+1)^{a^*+2},\nonumber\\
D(a)&=& D(r_++\frac{6}8(r-r_+),a),~~a>0,\nonumber\\
{\mathcal D}(a)&=&\{x: |{\rm Im}x|<a\},~~a>0,\nonumber\\
\Gamma(a)&=&\sum_{0<|k|\le K_+} |k|^{3n+(n+1)\tau}e^{-|k|\frac{a}8},
~~a>0,\nonumber\\
D_+&=& D(r_+,s_+),\nonumber\\
{\mathcal D}_+&=& {\mathcal D}(r_+)=\{x: |{\rm Im}x|<r_+\},\nonumber\\
D_{i}&=&D(r_++\frac{i-1}{8}(r-r_+),is_+),~~ i=1,2,\cdots,8,\nonumber
\end{eqnarray}
where $a^*$ is a constant such that $(\frac{10}9)^{a^*}>2.$

\subsection{Truncating perturbations}

Consider the  Taylor-Fourier series of $P$:
$$
P=\displaystyle \sum_{i\in Z_+^l,j\in Z_+^{2m}, k\in
Z^n}p_{kij}y^iz^je^{\sqrt{-1}\langle k,x\rangle}
$$
and consider the truncation
\begin{eqnarray}
R&=& \displaystyle \sum_{|i|+|j|\le 2, |k|\le
K_+}p_{kij}y^iz^je^{\sqrt{-1}\langle k,x\rangle}
=\displaystyle\sum_{|k|\le K_+}
         (P_{k00}+\langle P_{k10},y \rangle\nonumber\\
&&+\langle P_{k01},z \rangle+\langle y,P_{k20}y \rangle
 +\langle y,P_{k11}z \rangle+\langle z,P_{k02}z \rangle)
        e^{\sqrt{-1}\langle k,x \rangle}.\label{2.3}
\end{eqnarray}

\begin{lemma}
Assume that

{\bf H1)}
$$
\displaystyle\int^\infty_{K_+}\lambda^n e^{-\lambda
\frac{r-r_+}{8}}{\rm d}\lambda\leq \varepsilon.
$$
Then we have
$$
|\partial_\lambda^{l}(P-R)|_{D_8} \le \cdot\gamma^{n+1}
s^2\varepsilon^2,~~ |\partial_\lambda^{l}R|_{D_8} \le
\cdot\gamma^{n+1} s^2\varepsilon,~~~~|l|\le n.
$$
\end{lemma}

\noindent {\it Proof.} Let
\begin{eqnarray*}
I&=&\sum_{|k|>K_+}p_{kij}y^iz^j e^{\sqrt{-1}\langle k,x\rangle},\\
II&=&\sum_{|k|\le K_+,|i|+|j|> 2}p_{kij}y^iz^j e^{\sqrt{-1}\langle k,x\rangle}\\
&=&\int\frac{\partial^{(p,q)}} {\partial y^p\partial
z^q}\sum_{|k|\leq K_+,|i|+|j|> 2} p_{kij}e^{\sqrt{-1}\langle
k,x\rangle}y^i z^j{\rm d}y{\rm d}z,
\end{eqnarray*}
where $\displaystyle \int$ is the obvious anti-derivative of
$\displaystyle\frac{\partial^{(p,q)}} {\partial y^p\partial z^q}$
for $|p|+|q|=3$. Clearly,
$$
P-R=I+II.
$$

Since, by Cauchy's estimate,
$$
|\displaystyle \sum_{i\in Z_+^l,j\in Z_+^{2m}}\partial_\lambda^{l}
p_{kij}y^iz^j|\le|\partial_\lambda^{l} P|_{D(r,s)}e^{-|k|r} \le
\gamma^{n+1} s^2\varepsilon e^{-|k|r},~~|l|\le n,
$$
from H1) we get that
\begin{eqnarray}
|\partial_\lambda^{l} I|_{D_8}&\leq&\displaystyle\sum_{|k|> K_+}
\gamma^{n+1} s^2\varepsilon e^{-|k|r} e^{|k|(r_++\frac78(r-r_+))}\nonumber\\
&\leq& \gamma^{n+1} s^2\varepsilon\displaystyle\sum_{\kappa=
K_+}^{\infty}\kappa^n e^{-\kappa\frac{r-r_+}8} \leq \gamma^{n+1}
s^2\varepsilon\displaystyle\int^\infty_{K_+}\lambda^ne^{-\lambda\frac{r-r_+}{8}}
{\rm d}\lambda\nonumber\\
&\leq& \gamma^{n+1} s^{2}\varepsilon^2,~~~|l|\le n.\nonumber
\end{eqnarray}

It follows that
$$
|\partial_\lambda^{l}(P-I)|_{D_8}\le|\partial_\lambda^{l}
P|_{D(r,s)} +|\partial_\lambda^{l} I|_{D_8}\le \cdot\gamma^{n+1}
s^2\varepsilon,~~ |l|\le n.
$$

By  Cauchy's estimate we obtain
\begin{eqnarray}
|\partial_\lambda^{l} II|_{D_8}&\leq&
\left|\displaystyle\int\displaystyle\frac{\partial^{(p,q)}}
{\partial y^p\partial z^q}\displaystyle\sum_{|k|\leq K_+,|i|+|j|> 2}
\partial_\lambda^{l} p_{kij}
e^{\sqrt{-1}\langle k,x\rangle}y^i z^j{\rm d}y{\rm d}z\right|_{D_8}\nonumber\\
&\leq&
\left|\displaystyle\int\displaystyle\left|\frac{\partial^{(p,q)}}
{\partial y^p\partial z^q}\partial_\lambda^{l}(P-I-R)\right|_{D_{*}}
{\rm d}y{\rm d}z\right|_{D_8}\nonumber\\
&\le& \cdot\frac1{s^3}\gamma^{n+1}s^2\varepsilon\left|\int{\rm
d}y{\rm d}z\right|_{D_8} \leq
\cdot\frac1{s^3}\gamma^{n+1}s^2\varepsilon s_{+}^{3}
 \leq \cdot\gamma^{n+1}s^2\varepsilon^2,~~|l|\le n.\nonumber
\end{eqnarray}
Thus,
$$
|\partial_\lambda^{l} (P-R)|_{D_8}\leq c\gamma^{n+1}
s^2\varepsilon^2,
$$
and therefore,
$$
|\partial_\lambda^{l}R|_{D_8}\leq |\partial_\lambda^{l}(P-R)|_{D_8}+
|\partial_\lambda^{l}P|_{D_8}\leq\cdot \gamma^{n+1}
s^2\varepsilon,~~~|l|\leq n.
$$

\subsection{Transformation and  homogeneous equation}

Write ${\mathcal M}$ into blocks
\[
{\mathcal M}(x,\lambda)=\left(\begin{array}{cc}A & B\\ B^{\top} &
M\end{array}\right),
\]
where
$$
A(x, \lambda)=\sum_{k\in Z^n}A_ke^{\sqrt{-1}\langle k,x\rangle},~~
B(x,\lambda)=\sum_{k\in Z^n}B_ke^{\sqrt{-1}\langle  k,x\rangle},~~
M(x,\lambda)=\sum_{k\in Z^n}M_ke^{\sqrt{-1}\langle  k,x\rangle}
$$
are $l\times l$, $l\times 2m$, ${2m}\times {2m}$ minors of
${\mathcal M}$ respectively.

To transform (\ref{2.1}) into the Hamiltonian in the next KAM cycle,
we will construct the averaging transformation as the time 1-map
$\phi_F^1$ of the flow generated by a Hamiltonian $F$. To this end,
suppose $F$ has the following form:
\begin{equation}\label{2.4}
F=\sum_{0<|k|\le K_+}(f_{k0}+\langle f_{k1},y\rangle+\langle
F_{k1},z\rangle)e^{\sqrt{-1}\langle k,x\rangle}+\langle
F_{01},z\rangle.
\end{equation}

As in [7], to be able to keep the first $n_0$ components of the
toral pseudo-frequencies,  we shall also find a $Y_*\in R^{n_0}$ so
that the translation of coordinate
$$
 \phi:x\to x,~~~~y\to y+{Y_*\choose 0},~~~z\to z
$$
removes all possible drifts among the first $n_0$ components of the
new toral pseudo-frequencies.

We introduce the following notations:
\begin{eqnarray}
[A]&=&\left(\begin{array}{ll}U&D\\D^\top&V\end{array}\right),\nonumber\\
R'&=&\sum_{0<|k|\le K_+} (\langle y,P_{k20}y\rangle+\langle
y,P_{k11}z\rangle+ \langle z,P_{k02}z\rangle)e^{\sqrt{-1}\langle
  k,x\rangle} \nonumber\\
~~&&+[R] -\langle P_{001},z\rangle+\sum_{|k|\le K_+}
\langle B_{-k}JF_{k1}, y \rangle,\label{R'}\\
R_t&=&(1-t)\{N, F\}+ R,\label{Rt}\\
y_*&=&{Y_*\choose 0},\nonumber
\end{eqnarray}
where $U$, $D$, $V$ are the $n_0\times n_0$, $n_0\times (l-n_0)$,
$(l-n_0)\times (l-n_0)$ minors of $[A]$ respectively.

Denote
$$
\Phi_+=\phi^1_F\circ \phi.
$$

Then we have
\begin{eqnarray}
H_+&=&H\circ\Phi_+=H\circ\phi_F^1\circ\phi
=(N+R)\circ{\phi_F^1}\circ\phi+(P-R)\circ\phi_F^1\circ\phi\nonumber\\
&=& (N+R')\circ\phi-\langle y_*,(A-[A])y\rangle-\langle y_*,Bz\rangle\nonumber\\
&&+(\{N,F\}+R-R')\circ\phi+\langle y_*,(A-[A])y\rangle+\langle
y_*,Bz\rangle-Q
\nonumber\\
&&+\displaystyle\int^1_0\{R_t,F\}\circ\phi^t_F \circ\phi{\rm d}t
+(P-R)\circ\phi_F^1\circ\phi+Q,\nonumber
\end{eqnarray}
where $Q$ is to be determined in the following.

As in [7], we need to choose a function $Q$ such that both equations
\begin{eqnarray}
&&(\{N,F\}+R-R')\circ\phi-Q+\langle y_*,(A-[A])y\rangle+\langle
y_*,Bz\rangle=0,
\label{2.5}\\
&& {\rm diag}(U, O)y_*={\rm diag}(I_{n_0}, O)(-P_{010}-\sum_{|j|\le
K_+}B_{-j}JF_{j1})\label{yy}
\end{eqnarray}
are solvable.  If this is the case, we then arrive at that
\begin{eqnarray}
H_+&=&N_++P_+,\nonumber\\
N_+&=&e_++\langle\Omega_+(\lambda),y\rangle+ \frac12\left\langle
{y\choose z},{\mathcal M}^+{y\choose z}\right\rangle+h_0(x,y,z,
\lambda)\nonumber\\
~~~~&=&e_++\langle\Omega_+(\lambda),y\rangle+ \frac12\left\langle
{y\choose z}, \left(\begin{array}{cc}A^+ & B^+\\ {B^+}^{\top} & M^+
\end{array}\right){y\choose z}\right\rangle+h_0(x,y,z,\lambda),\nonumber
\end{eqnarray}
where
\begin{eqnarray}
&&e_+=e+P_{000}+\langle \Omega, y_*\rangle+\frac12\langle y_*,[A]y_*
\rangle,
\label{2.01}\\
&&\Omega_+=\Omega+{\rm diag}(O, I_{n-n_0})([A]y_*+P_{010}+
\sum_{|k|\le K_+}
B_{-k}JF_{k1}),\label{2.02}\\
&&\omega_+=-E^\top\Omega_+,\\
&&A^+=A+\sum_{|k|\le K_+}2P_{k20}e^{\sqrt{-1}\langle k,x\rangle},\label{2.02a}\\
&&B^+=B+\sum_{|k|\le K_+}P_{k11} e^{\sqrt{-1}\langle k,x\rangle},\label{2.02b} \\
&&M^+=M+\sum_{|k|\le K_+}2P_{k02} e^{\sqrt{-1}\langle k,x\rangle},\label{2.02c}\\
&&P_+=\displaystyle\int^1_0\{R_t,F\}\circ \phi_F^t\circ\phi {\rm
  d}t+(P-R)\circ\phi_F^1\circ\phi\nonumber\\
&&\;~~~~~+\frac 12\langle y_*, (A-[A])y_*\rangle
  +h_0(x,y+y_*,z,\lambda)-h_0(x,y,z,\lambda)\nonumber\\
&&\;~~~~~+\sum_{|k|\le K_+}(\langle y_*,P_{k20}y_*\rangle
  +\langle y_*,2P_{k20}y\rangle+\langle y_*, P_{k11}z\rangle)
  e^{\sqrt{-1}\langle k,x\rangle}+Q.\label{2.6}
\end{eqnarray}

We now consider the equations (\ref{2.5}) and (\ref{yy}). By careful
observation of (\ref{2.5}), we suppose that $Q$ has the following
form:
\begin{align}
Q=&(\sum_{0<|k|\le K_+} (-\langle\frac 12\partial_x\langle
y,A(x,\lambda)y\rangle +\partial_x\langle
y,B(x,\lambda)z\rangle\nonumber\\
&+\frac 12\partial_x\langle z,M(x,\lambda)z\rangle+\partial_x
h_0(x,y,z,\lambda), E^\top f_{k1}\rangle\nonumber\\
&+\sqrt{-1}\langle E k,A(x,\lambda)y+B(x,\lambda)z\nonumber\\
&+\partial_y h_0(x,y,z,\lambda)\rangle(f_{k0}+\langle
f_{k1},y\rangle+\langle F_{k1},z\rangle)\nonumber\\
&+\sqrt{-1}\langle\frac 12\partial_x\langle y,A(x,\lambda)y\rangle
+\partial_x\langle y,B(x,\lambda)z
\rangle\nonumber\\
&+\frac 12\partial_x\langle z,M(x,\lambda)z\rangle+\partial_x
h_0(x,y,z,\lambda),Ck\rangle(f_{k0}+\langle
f_{k1},y\rangle+\langle F_{k1},z\rangle)) e^{\sqrt{-1}\langle k,x\rangle}\nonumber\\
&+\sum_{|k|> K_+}(\langle B_kJF_{01}, y\rangle +\langle
M_kJF_{01},z\rangle) e^{\sqrt{-1}\langle k,x\rangle}\nonumber\\
&+\sum_{|k|> K_+, 0<|j|\le K_+}(\langle
B_{k-j}JF_{j1},y\rangle+\langle
M_{k-j}JF_{j1},z\rangle)  e^{\sqrt{-1}\langle k,x\rangle}\nonumber\\
&+\sum_{ 0<|k|\le K_+}\langle \partial_z
h_0(x,y,z,\lambda)JF_{k1},z\rangle
e^{\sqrt{-1}\langle k,x\rangle}) \circ\phi\nonumber\\
&+\sum_{0<|k|\le K_+}(\sqrt{-1}\langle
k,E^\top\Omega(\lambda)\rangle\langle f_{k1},y_*\rangle +\langle
P_{k10},y_*\rangle)e^{\sqrt{-1}\langle k,x\rangle}\nonumber\\
&+\sum_{|k|> K_+}(\langle y_*,A_ky\rangle +\langle y_*,B_k
z\rangle)e^{\sqrt{-1}\langle k,x\rangle}\nonumber\\
&+\sum_{0<|k|\leq K_+,0\leq |j|\leq K_+}\langle
y_*,B_{k-j}JF_{j1}\rangle e^{\sqrt{-1}\langle k,x\rangle}. \label{Q}
\end{align}

Substituting (\ref{2.3})--(\ref{R'}) and (\ref{Q}) into (\ref{2.5})
and comparing coefficients, from equations (\ref{2.5}) and
(\ref{yy}) we obtain the following linear equations for all
$0<|k|\le K_+$:
\begin{eqnarray}
& &\sqrt{-1}\langle k,\omega(\lambda)\rangle f_{k0}=P_{k00},\label{2.001}\\
& &\sqrt{-1}\langle k,\omega(\lambda)\rangle
f_{k1}=P_{k10}+A_ky_*+\sum_{|j|
\le K_+}B_{k-j}JF_{j1},\label{2.002}\\
& &\sqrt{-1}\langle k,\omega(\lambda)\rangle F_{k1}-[M]J F_{k1}
=\sum_{0<|j|\le K_+, j\ne k}M_{k-j}JF_{j1}\nonumber\\
&&\qquad\qquad\qquad\qquad\qquad\qquad\qquad+P_{k01}+B_k^{\top}y_*+M_kJF_{01},
\label{2.003}\\
& &[M]J F_{01}=-P_{001}-\sum_{0<|j|\le
K_+}M_{-j}JF_{j1}-[B]^{\top}y_*,
\label{2.004}\\
&& {\rm diag}(U, O)y_* ={\rm
diag}(I_{n_0},O)(-P_{010}-\sum_{0<|j|\le
K_+}B_{-j}JF_{j1}-[B]JF_{01}), \label{y}
\end{eqnarray}
where $\omega(\lambda)=-E^\top(\lambda)\Omega(\lambda).$

Denote
\begin{equation}\label{2.10} {\Lambda}_+=\{\lambda\in
{\Lambda}: |\langle
k,\omega(\lambda)\rangle|>\frac{\gamma}{|k|^\tau}, ~~0<|k|\le K_+\}.
\end{equation}

If we assume that

{\bf H2)}
$$
\displaystyle|\partial^{l}_\lambda (\mathcal M-\mathcal
M^0)|_{{\mathcal D}(r)\times \Lambda} \le \varepsilon_0^{\frac14},
~~~|l|\le n,
$$
then as in [7], (\ref{2.001})--(\ref{y}) can be equivalently written
into the following system form:
\begin{equation}\label{2.9}
(\Lambda-\mathcal A){\mathcal F}={\mathcal P},
\end{equation}
where $\Lambda$ and $\mathcal A$ are defined as in [7].

\subsection{Estimate on $(\Lambda-{\mathcal A})^{-1}$}
As in [7], by the hyperbolicity of $J[M^0]$ and the definition of
$\eta$, we can prove that
\begin{equation}\label{Lambda0}
|(\Lambda^0)^{-1}|_{{\Lambda}_0}\le \frac{\rho_0}2,~~~ |{\mathcal
A}^0|_{{\Lambda}_0}\le \frac1{\rho_0},~~~ |(\Lambda^0-{\mathcal
A}^0)^{-1}|_{{\Lambda}_0}<2\rho_0.
\end{equation}

\begin{lemma} Assume  {\rm H2)} and also that

{\rm \bf H3)}
$$
|\partial^{l}_\lambda {\mathcal A}-\partial^{l}_\lambda{\mathcal
A}^0|_{\Lambda}<\varepsilon_0^{\frac14}.
$$
Then for $\varepsilon_0$  sufficiently small, $\mathcal
L=\Lambda-{\mathcal A}$ is non-singular on $\Lambda$, and moreover,
the following holds:
\[
|\partial^{l}_\lambda\mathcal L^{-1}|_{\Lambda}\le \cdot
K_+^n,~~~~|l|\le n.
\]
\end{lemma}

\noindent {\it Proof.} Similar to Lemma 3.2 of [7], we have $
|\partial_\lambda\mathcal L^{-1}|_{\Lambda}\le \cdot K_+. $ By
induction,
$$
|\partial_\lambda^{l}\mathcal L^{-1}|_{\Lambda}\le \cdot
K_+^n,~~~~|l|\le n.
$$

Above all, by the hypotheses H2) and H3), the linear system
(\ref{2.9}) can be uniquely solved on $\Lambda_+$ to yield  smooth
functions $f_{k0}, f_{k1}, F_{k1}, F_{01}, y_*$, $0<|k|\le K_+$.

\subsection{Estimates on the transformation}

Denote
$$
\zeta=K_+^{n+2}\Gamma(r-r_+)^2.
$$

\begin{lemma} Assume {\rm H2)}.
Then the following holds for all $|l|\le n$:

{\rm 1)}
 $|\partial_\lambda^{l}y_*|_{\Lambda_+}\le \cdot \gamma^{n+1} s\varepsilon \zeta$;

{\rm 2)} On $D(s)\times {\Lambda}_+$,
$$
|\partial_\lambda^{l}F|,~
|\partial_\lambda^{l}F_x|,~s|\partial_\lambda^{l}F_y|,~
s|\partial_\lambda^{l}F_z|\le \cdot s^2\varepsilon \zeta;
$$

{\rm 3)} On $D(s)\times {\Lambda}_+$,
$$
|\partial_\lambda^{l}D^iF|\le \cdot\varepsilon \zeta,~~|i|\ge 2.
$$
\end{lemma}

\noindent {\it Proof.} The proof is similar to that in [7].

\begin{lemma} Assume {\rm H2), H3)} and also that

{\bf H4)}
$$
s\varepsilon \zeta<\frac{1}{8}(r-r_+),~~~ s\varepsilon\zeta < s_+.
$$
Let $\phi^t_F$  be the flow generated by $F$. Then the following
holds:

{\rm 1)}  For all $0\leq t\leq 1$, $\phi^t_F:D_2\rightarrow D_3,~
\phi: D_1 \rightarrow D_2$ are well defined, real analytic and
depend smoothly on $\lambda\in \Lambda_+$, i.e.,
$\Phi_+=\phi_F^1\circ \phi:D_+\rightarrow D$;

{\rm 2)} $|\partial_\lambda^{l}(\phi_F^t-id)|_{{D(s)}\times
\Lambda_+}\le \cdot s\varepsilon \zeta,~
|\partial_\lambda^{l}D^i(\Phi_+-id)|_{ {\tilde D}_+\times \Lambda_+}
\leq \cdot\varepsilon \zeta,$ for all $|l|\le n,~i\ge 0,~0\le t\le
1$, where $D=\partial_{(x,y,z)}$.
\end{lemma}

\noindent {\it Proof.} Let $\lambda\in {\Lambda}_+$.

1) It is easy to see that $\phi:D_1\rightarrow D_2$ holds by Lemma
2.3 1) and H4).

We note that
\begin{equation}\label{phit}
\phi^t_F=   {\rm id}+\int^t_0X_F\circ\phi^\xi_F{\rm d}\xi,
\end{equation}
where
$$
X_F=\tilde{I}(\lambda)\nabla F=(E(\lambda)F_x, -E^\top(\lambda)
F_y+C(\lambda)F_x, JF_z)^{\top}.
$$
Denote $\phi^t_{F1},\phi^t_{F2}, \phi^t_{F3}$ as components of
$\phi_F^t$ in $y,x,z$ planes respectively.  For any $(x,y,z)\in
D_2$,  let $t_*={\rm sup} \{t\in [0,1]: \phi_F^t(x,y,z)\in D_3$. By
making $\varepsilon_0$ small, we have that $D_3\subset D(s)$. It
follows from H4) and Lemma 2.3 that
\begin{eqnarray}
|\phi_{F1}^t(x,y,z)| &\le& |y|+|\int_0^t
       E(\lambda)F_x\circ\phi_F^{\xi}{\rm d}\xi|
    \leq |y|+\cdot|F_x|_{D(s)}
    \leq 2s_+ +\cdot s^2\varepsilon\zeta<3s_+,\nonumber\\
|\phi_{F2}^t(x,y,z)| &\le& |x|+ |\int_0^t
      (-E^\top(\lambda)F_y+C(\lambda)F_x)\circ\phi_F^{\xi}{\rm d}\xi|
      \leq|x|+\cdot(|F_x|+|F_y|)_{D(s)}\nonumber\\
      &\leq& r_++\frac18(r-r_+)+ \cdot s\varepsilon \zeta\nonumber\\
    &<& r_++\frac28(r-r_+),\nonumber\\
|\phi_{F3}^t(x,y,z)| &\le& |z|+|\int_0^t
       JF_z\circ\phi_F^{\xi}{\rm d}\xi|
       \leq |z|+|F_z|_{D(s)}
       \leq 2s_++\cdot s\varepsilon \zeta
       < 3s_+,      \nonumber
\end{eqnarray}
i.e., $\phi_F^t(x,y,z)\in D_3$ for all $0\le t\le t_*$. Thus,
$t_*=1$ and 1) holds.

2) By Lemma 2.3 and (\ref{phit}), we immediately have
$$
|\phi_F^t-id|_{{D(s)}}\le \cdot s\varepsilon \zeta.
$$

Differentiating  (\ref{phit}) with respect to $\lambda$ yields
\begin{align*}
\partial_\lambda\phi^t_F=&
\int^t_0X_F\circ\phi^\xi_F\partial_\lambda\phi_F^\xi{\rm d}\xi
+\int^t_0 (\partial_\lambda X_F)\circ \phi_F^\xi{\rm d}\xi\\
=&\int^t_0(E(\lambda)F_x, -E^\top(\lambda) F_y+C(\lambda)F_x,
JF_z)^{\top}\circ\phi^\xi_F\partial_\lambda\phi_F^\xi{\rm d}\xi\\
&+\int^t_0 \partial_\lambda(E(\lambda)F_x, -E^\top(\lambda)
F_y+C(\lambda)F_x, JF_z)^{\top} \circ\phi^\xi_F{\rm d}\xi.
\end{align*}
It follows from Lemma 2.3 and Gronwall's inequality that
$$
|\partial_\lambda\phi^t_F|_{D(s)}\le \cdot s\varepsilon \zeta.
$$
By induction, we have
$$
|\partial_\lambda^{l}\phi^t_F|_{D(s)}\le \cdot s\varepsilon
\zeta,~~|l|\le n.
$$

The estimates for $\Phi_+$ follow from a similar application of
Lemma 2.3 and  Gronwall's inequality, and the identity
$$
\Phi_+-id=(\phi^1_F-id)\circ \phi+\left(\begin{array}{lll}
0\\y_*\\0\end{array}\right).
$$
We omit the details.

\subsection{Estimate on $N_+$}

We first estimate the new normal form.

\begin{lemma} For the new normal form, we have
  the following holds for all $|l|\le n$:
\begin{eqnarray}
|\partial_\lambda^{l}(e_+-e)|_{{\Lambda}_+}&\le&
\cdot~\gamma^{n+1}s\varepsilon\zeta,\nonumber\\
|\partial_\lambda^{l}(\Omega_+-\Omega)|_{{\Lambda}_+}&\le&
\cdot~\gamma^{n+1} s\varepsilon \zeta,  \nonumber\\
|\partial_\lambda^{l}(\omega_+-\omega)|_{{\Lambda}_+}&\le&
\cdot~\gamma^{n+1} s\varepsilon \zeta,  \nonumber\\
|\partial_\lambda^{l}({\mathcal M}^+-{\mathcal M})|_{{\mathcal
D}_+\times {\Lambda}_+}&\le&
\cdot~\gamma^{n+1}\varepsilon\Gamma(r-r_+). \nonumber
\end{eqnarray}
\end{lemma}

\noindent {\it Proof.} First, by  Cauchy's estimate we have
\begin{eqnarray}\label{pij}
|\partial_\lambda^{l} P_{kij}|_{\mathcal O}&\le& \cdot s^{-(i+j)}
|\partial_\lambda^{l} P|_{D(r,s)\times\mathcal O}e^{-|k|r}\nonumber\\
&\le& \cdot\gamma^{n+1} s^{2-i-j}\varepsilon e^{-|k|r},~~~|k|\ge
0,~i,j=0,1,2.
\end{eqnarray}
Then from (\ref{2.01})--(\ref{2.02c}) and (\ref{pij}) the Lemma
immediately follows.

\subsection{Frequency property}

\begin{lemma} Assume that

{\bf H5)}
$$
\gamma^{n+1}s\varepsilon\zeta K_+^{\tau+1}<\gamma-\gamma_+.
$$
Then
$$
|\langle k,\omega_+(\lambda)\rangle|>\frac{\gamma_+}{|k|^{\tau}},
$$
for all $\lambda\in {\Lambda}_+$ and $0<|k|\le K_+$.
\end{lemma}

\noindent {\it Proof.} By H5) and Lemma 2.5, one has
\begin{eqnarray}
|\langle k,\omega_+(\lambda)\rangle|&=&|\langle
k,\omega(\lambda)\rangle+\langle
k,\omega_+(\lambda)-\omega(\lambda)\rangle|\nonumber\\
&\ge&|\langle k,\omega(\lambda)|-\gamma^{n+1}s\varepsilon \zeta K_+\nonumber\\
&\ge&\frac{\gamma}{|k|^\tau}-\frac{\gamma-\gamma_+}{|k|^\tau}=
\frac{\gamma_+}{|k|^\tau},
\end{eqnarray}
as desired.

\subsection{Estimate on the new perturbation}

Denote
\begin{equation}\label{Delta}
\Delta=\cdot s^3\varepsilon^2\zeta^2+\cdot\gamma^{n+1}
s^{2}\varepsilon^2 \zeta^2+\cdot s_+s^2\varepsilon\zeta.
\end{equation}

\begin{lemma} Assume {\rm H1)--H4)}. Then
$ |\partial_\lambda^{l}P_+|_{D_+}\le \Delta,~|l|\le n. $ Thus, if

{\rm \bf H6)}
$$
\Delta\le \gamma_+^{n+1} s_+^2\varepsilon_+,
$$
then
\begin{equation}\label{P1}
|\partial_\lambda^{l}P_+|_{D_+}\le \gamma_+^{n+1}
s_+^2\varepsilon_+.
\end{equation}
\end{lemma}

\noindent {\it Proof.} Let $|l|\le n$, $\lambda\in \Lambda_+$. By
(\ref{2.6}), we have that
\begin{equation}\label{P}
P_+=W_0\circ\phi+W_1+Q+q+(P-R)\circ\Phi_+,
\end{equation}
 where
\begin{eqnarray*}
W_0&=&\int^1_0\{R_t,F\}\circ \phi_F^t{\rm d}t,\\
W_1&=&\frac12\langle y_*,(A-[A])y_*\rangle+\sum_{|k|\le K_+}(\langle
y_*,P_{k20}y_*\rangle
  +\langle y_*,2P_{k20}y\rangle+\langle y_*, P_{k11}z\rangle)
  e^{\sqrt{-1}\langle k,x\rangle},\\
q&=&h_0(x,y+y_*,z,\lambda)-h_0(x,y,z,\lambda).
\end{eqnarray*}

1) We first estimate $(P-R)\circ\Phi_+$.

By Lemma 2.4 1) and Lemma 2.1, we have
\begin{equation}\label{P-R}
|\partial_\lambda^l(P-R)\circ\Phi_+|_{D_+}\leq
|\partial_\lambda^l(P-R)|_{D_3} \le
\cdot\gamma^{n+1}s^2\varepsilon^2.
\end{equation}

2) Then we give the estimate of $q$.

Following the Taylor series expansion, H4) and Lemma 2.3 1), we
obtain
\begin{align}
|\partial_\lambda^l
q|_{D_+}&=|\partial_\lambda^l(h'_0(y)y_*)+\frac1{2!}y_*h_0^{(2)}y_*+
\frac1{3!}h_0^{(3)}y_*^3+o(y_*^3)|_{D_+}\nonumber\\
&\le s_+^2|y_*|+s_+|y_*|^2+|y_*|^3 \le \cdot s_+^2|y_*|\le \cdot
\gamma^{n+1}s_+^2s\varepsilon\zeta.\label{q}
\end{align}

3) Then we estimate $W_1$.

By Lemma 2.3 1), (\ref{2.6}) and H4), we have that
\begin{eqnarray}
|\partial_\lambda^l W_1|_{D_+}&\le&\cdot|y_*|^2+\sum_{|k|\le K_+}
(|y_*|^2\gamma^{n+1}\varepsilon+s_+|y_*|\gamma^{n+1}\varepsilon)e^{-|k|\frac{r-r_+}2}
\nonumber\\
&\le&\cdot|y_*|^2+\sum_{|k|\le K_+}\cdot
s_+|y_*|\gamma^{n+1}\varepsilon e^{-|k|\frac{r-r_+}2}\nonumber\\
&\le &\cdot\gamma^{n+1}s^2\varepsilon^2\zeta^2+
\cdot s_+\gamma^{n+1}s\varepsilon\zeta\gamma^{n+1}\varepsilon\Gamma\nonumber\\
&\le& \cdot\gamma^{n+1}s^2\varepsilon^2\zeta^2.\label{W1}
\end{eqnarray}

4) Next, we give the estimate of $Q$.

By a similar computation to [7], and noting that
$|E(\lambda)|,|C(\lambda)|\le c$ for some constant c, we obtain that
\begin{equation}\label{QE}
|\partial_{\lambda}^{l}Q|_{D_+} \le \cdot
s_+s^2\varepsilon\zeta+\cdot\gamma^{n+1}s^2\varepsilon^2\zeta^2.
\end{equation}

5) Now we can estimate $W_0\circ\phi$.

We can obtain  the estimate of $W_0\circ\phi$ as in [7]:
\begin{equation*}
|\partial_\lambda^{l}W_0\circ\phi|_{D_+}\le \cdot
s^3\varepsilon^2\zeta^2+\cdot\gamma^{n+1} s^{2}\varepsilon^3
\zeta^3+\cdot\gamma^{n+1} s^{2}\varepsilon^2 \zeta^2.
\end{equation*}

It will be proved that $\varepsilon\zeta\le 1$ later, so we have
\begin{equation}\label{W0}
|\partial_\lambda^{l}W_0\circ\phi|_{D_+}\le \cdot
s^3\varepsilon^2\zeta^2+\cdot\gamma^{n+1} s^{2}\varepsilon^2
\zeta^2.
\end{equation}

Above all, it follows from (\ref{P-R}), (\ref{q}), (\ref{QE}),
(\ref{P}), (\ref{W0}), (\ref{W1}) that
$$
|\partial_\lambda^{l}P_+|_{D_+}\le \cdot
s^3\varepsilon^2\zeta^2+\cdot\gamma^{n+1} s^{2}\varepsilon^2
\zeta^2+\cdot s_+s^2\varepsilon\zeta.
$$
So by H6), (\ref{P1}) holds. This completes the proof of the Lemma.

This completes one cycle of KAM steps.

\section{Iteration Lemma}
\setcounter{equation}{0}

Consider (\ref{1.1}) and let $r_0,s_0, \varepsilon_0, \gamma_0,
{\Lambda}_0, H_0, N_0, e_0, \Omega_0, {{\mathcal M}}^0, A^0, B^0,
M^0, {{\mathcal A}}^0, h_0, P_0$ be defined in Section~2 and let
$D_0=D(r_0, s_0)$, ${\mathcal D}_0=\{x: |{\rm Im}x|<r_0\}$, $K_0=0$,
$\Phi_0=id$. We define the following sequences inductively for all
$\nu=1,2,\cdots:$
\begin{eqnarray}
H_\nu&=&H_\nu(x,y,z,\lambda)=N_{\nu}+P_{\nu},\nonumber\\
N_{\nu}&=&e_{\nu}+\langle
\Omega_{\nu},y\rangle+\frac12\left\langle{y\choose z},
{\mathcal M}^{\nu}{y\choose z}\right\rangle+h_0(x,y,z,\lambda),\nonumber\\
{\mathcal M}^{\nu}&=&\left(\begin{array}{cc}A^{\nu} & B^{\nu}\\
{(B^{\nu})}^{\top} & M^{\nu}
\end{array}\right),\nonumber\\
\varepsilon_{\nu}&=&\varepsilon_{\nu-1}^{\frac{10}9},\nonumber\\
r_\nu&=   &r_0\left(1-\sum_{i=1}^\nu\frac1{2^{i+1}}\right),\nonumber\\
s_\nu&=   &\frac18\alpha
s_{\nu-1},~~\alpha_{\nu-1}=\varepsilon_{\nu-1}^{\frac13},
\nonumber\\
\gamma_\nu&=&   \gamma_0\left(1-\sum_{i=1}^\nu\frac1{2^{i+1}}\right),\nonumber\\
K_{\nu}&=
&\left(\left[\log\frac{1}{s_{\nu-1}}\right]+1\right)^3,~~\nu\ge 1,
\nonumber\\
{\Lambda}_{\nu}&=&\{\lambda\in {\Lambda}_{\nu-1}:|\langle
k,\Omega_{\nu-1}(\lambda)\rangle|>\frac{\gamma_{\nu-1}}{|k|^\tau},
0<|k|\le
K_{\nu}\},~\nu\ge 1,\nonumber\\  D_{\nu}&=   &D(r_{\nu},s_{\nu}),\nonumber\\
\mathcal D_{\nu}&=& \{x: |{\rm Im}x|< r_{\nu}\}.\nonumber
\end{eqnarray}

\begin{lemma}\hspace{-2mm}{\rm (Iteration Lemma)} \ If
$\varepsilon_0=\varepsilon_0(r_*,\sigma_0,U^0)$  is sufficiently
small, then the following holds for all $|l|\le n;\nu=1,2,\cdots$.

{\rm 1)}  There is a transformation $\Phi_{\nu}:{D}_{\nu}\times
\Lambda_{\nu} \longrightarrow {D}_{\nu-1}$, which  is  symplectic
and analytic in $(x,y,z)\in D_{\nu+1}$, and smooth in $\lambda\in
\Lambda_{\nu+1}$, such that $
 H_{\nu}=H_{\nu-1}\circ\Phi_{\nu}
$ and
\begin{equation}\label{310}
 |\partial_\lambda^{l}D^i(\Phi_{\nu}-id)|_{D_{\nu}\times \Lambda_{\nu}}
 \le \cdot\varepsilon_{\nu-1}\zeta_{\nu-1},\ i\ge 0.
 \end{equation}

{\rm 2)}
  \begin{eqnarray}
 &&|\partial_\lambda^{l}(e_{\nu}-e_0)|_{\Lambda_{\nu}}\le
\cdot\gamma_0^{n+1}s_0\varepsilon_0\zeta_0,\label{31}\\
 && |\partial_\lambda^{l}(e_{\nu}-e_{\nu-1})|_{\Lambda_{\nu}}
 \le \cdot\gamma_{\nu-1}^{n+1}s_{\nu-1}\varepsilon_{\nu-1}\zeta_{\nu-1},\label{32}\\
 &&  |\partial_\lambda^{l}(\Omega_{\nu}-\Omega_0)|_{\Lambda_{\nu}}\le
 \cdot\gamma_0^{n+1}s_0\varepsilon_0\zeta_0,\label{33}\\
 && |\partial_\lambda^{l}(\Omega_{\nu}-\Omega_{\nu-1})|_{\Lambda_{\nu}}\le
 \cdot\gamma_{\nu-1}^{n+1}s_{\nu-1}\varepsilon_{\nu-1}\zeta_{\nu-1},\label{34}\\
 &&  |\partial_\lambda^{l}(\omega_{\nu}-\omega_0)|_{\Lambda_{\nu}}\le
 \cdot\gamma_0^{n+1}s_0\varepsilon_0\zeta_0,\label{35}\\
 && |\partial_\lambda^{l}(\omega_{\nu}-\omega_{\nu-1})|_{\Lambda_{\nu}}\le
 \cdot\gamma_{\nu-1}^{n+1}s_{\nu-1}\varepsilon_{\nu-1}\zeta_{\nu-1},\label{36}\\
 &&  |\partial_\lambda^{l}({\mathcal M}^{\nu}-{\mathcal M}^0)|_{\mathcal D_{\nu}\times
 \Lambda_{\nu}}\le \cdot\gamma_0^{n+1}\varepsilon_0\zeta_0,\label{37}\\
&&|\partial_\lambda^{l}({\mathcal M}^{\nu}-{\mathcal
M}^{\nu-1})|_{\mathcal D_{\nu}\times\Lambda_{\nu}}\le
\cdot\gamma_{\nu-1}^{n+1}\varepsilon_{\nu-1}\zeta_{\nu-1},
\label{38}\\
&&|\partial_\lambda^{l}P_{\nu}|_{D_{\nu}\times\Lambda_{\nu}}\leq
\gamma_\nu^{n+1} s_{\nu}^2\varepsilon_{\nu}.\label{39}
\end{eqnarray}

{\rm 3)}
$(\Omega_{\nu}(\lambda))_{i}=\Omega_{i}(\lambda),~~i=1,2,\cdots,n_0$.

{\rm 4)}
\[
{\Lambda}_{\nu+1}=\{\lambda\in {\Lambda}_{\nu}: |\langle
k,\Omega_{\nu}(\lambda)\rangle|>\frac{\gamma_{\nu}}{|k|^\tau},
~K_{\nu}<|k|\le K_{\nu+1}\}.
\]
\end{lemma}

\noindent {\it Proof.} The proof amounts to the verification of
H1)--H6) for all $\nu$. The lemma will be proved by induction.

By definitions of $\varepsilon_\nu, s_\nu$, we clearly have
 \begin{eqnarray}
 \varepsilon_\nu&=&\varepsilon_0^{(\frac{10}9)^\nu},\label{3.04}\\
s_\nu&=&\left(\frac18\right)^\nu\varepsilon_0^{3((\frac{10}9)^\nu-1)}s_0.\label{3.05}
\end{eqnarray}

By making $\varepsilon_0$ small, we obtain that
\begin{eqnarray}
&&\quad\log (n+1)!+n(a^*+2)\log\left(\left[\log\frac
1{\varepsilon_\nu}\right]+1
\right)-\frac1{2^{\nu+5}}\left(\left[\log\frac
1{\varepsilon_\nu}\right]
+1\right)^{a^*+2}r_0\nonumber\\
&&\quad+(n+1)((\nu+5)\log2-\log r_0)\nonumber\\
&&\le \log (n+1)!+ n(a^*+2)\log\left(\log\frac
1{\varepsilon_\nu}+2\right)
-\left(\log\frac 1{\varepsilon_\nu}\right)^2r_0\nonumber\\
&&\quad+(n+1)((\nu+5)\log2-\log r_0)\nonumber\\
&&\le -\log\frac1{\varepsilon_\nu},\nonumber
\end{eqnarray}
where the first `$\le$' holds because of the choice of $a^*$
satisfying
$\frac1{2^{\nu+5}}\left(\log\frac1{\varepsilon_\nu}\right)^{a^*}\ge
1.$ Hence,
$$
\displaystyle\int_{K_{\nu+1}}^{\infty}\lambda^n e^{ -\lambda
\frac{r-r+}8} {\rm d}\lambda \leq
(n+1)!K_{\nu+1}^n\left(\frac{2^{\nu+5}}{r_0}
\right)^{n+1}e^{-K_{\nu+1}\frac{r_0}{2^{\nu+5}}}\le s_\nu.
$$
This verifies H1).

The verification of H2)--H3) is similar to that in [7], the reader
can refer to [7] for details.

To verify H4)--H6), we will prove $\varepsilon^{\frac19}\zeta^2\le
1$ at first. By making $\varepsilon_0$ sufficiently small, it
follows that
\begin{eqnarray}
\varepsilon^{\frac19}\zeta^2&=&\varepsilon_0^{\frac19(\frac{10}9)^\nu}
(K_+^{n+2}\Gamma(r-r_+)^2)^2\nonumber\\
&\le&\varepsilon_0^{\frac19(\frac{10}9)^\nu}\left(\log\frac1{\varepsilon}+1\right)
^{2(n+2)(a^*+2)}((3n+(n+1)\tau+1)!)^4\left(\frac{2^{\nu+5}}{r_0}\right)^{4(3n+(n+1)
\tau+1)}\nonumber\\
&\le
&\cdot\varepsilon_0^{\frac19(\frac{10}9)^\nu}\left(\log\frac1{\varepsilon}+1
\right)^{2(n+2)(a^*+2)}2^{4\nu(3n+(n+1)\tau+1)}\nonumber\\
&\le
&\cdot[\varepsilon_0^{\frac1{18}(\frac{10}9)^\nu}\left(\log\frac1{\varepsilon}
+1\right)^{2(n+2)(a^*+2)}][\varepsilon_0^{\frac1{18}{(\frac{10}9)^\nu}}2^{4\nu(3n
+(n+1)\tau+1)}]\nonumber\\
&\le & 1.\label{gamma}
\end{eqnarray}

Now, we can prove H4)--H6) easily. When $\varepsilon$ is
sufficiently small, the following hold:
\begin{align*}
&s\varepsilon\zeta\le
s\varepsilon^{\frac89}=\left(\frac18\right)^\nu\varepsilon_0^{3((\frac{10}9)^\nu-1)}
s_0\varepsilon_0^{\frac89(\frac{10}9)^\nu}<\left(\frac12\right)^{\nu+5}r_0=
\frac{r-r_+}8,\\
&s\varepsilon\zeta \le
s\varepsilon^{\frac89}<\frac18\varepsilon^{\frac13}s=s_+,\\
&\gamma^{n+1}s\varepsilon\zeta K_+^{\tau+1}\le^{(*)}
\gamma_0^{n+1}s\varepsilon^{\frac12}
<\frac{\gamma_0}{2^{\nu+2}}=\gamma-\gamma_+,
\end{align*}
i.e., H4, H5) hold, where $(*)$ holds similar to (\ref{gamma}).

At last, we give the proof of H6). By the smallness of
$\varepsilon_0$ and the choice of
$s_0=\left(\displaystyle\frac{\gamma_0}2\right)^{n+1}
\varepsilon_0^{\frac59}$, we have the following estimates:
\begin{align*}
&s^3\varepsilon^2\zeta^2=8^2\left(\frac18\right)^2\alpha^2s^2(\varepsilon^{\frac19}
\zeta^2)s\varepsilon^{\frac19}\varepsilon^{\frac{10}9}\le
8^2s\varepsilon^{\frac19}s_+^2\varepsilon_+\le
\gamma_+^{n+1}s_+^2\varepsilon_+,\\
&\gamma^{n+1}s^2\varepsilon^2\zeta^2\le \gamma^{n+1}s^2\varepsilon
\varepsilon^{\frac89}=8^2\gamma^{n+1}\varepsilon^{\frac19}s_+^2\varepsilon_+
\le \gamma_+^{n+1}s_+^2\varepsilon_+,\\
&s_+s^2\varepsilon\zeta\le
s_+s^2\varepsilon^{\frac89}=8s_+^2\varepsilon_+s\varepsilon^{-\frac59}\\
&\qquad\quad=8s_+^2\varepsilon_+\left(\frac18\right)^\nu\varepsilon_0^{3((
\frac{10}9)^\nu-1)}(\frac{\gamma_0}2)^{n+1}\varepsilon_0^{\frac59}\varepsilon_0^{
-\frac59(\frac{10}9)^\nu}\\
&\qquad\quad\le \gamma_+^{n+1}s_+^2\varepsilon_+.
\end{align*}
This verifies H6).

Above all, H1)--H6) hold for all $\nu=0,1,\cdots$, i.e., the KAM
step described in Section~2 is
 valid for all $\nu=0,1,\cdots$. Now, (\ref{32}), (\ref{34}),
(\ref{36}) and (\ref{38}) follow from Lemma 2.5; (\ref{31}),
(\ref{33}), (\ref{35}) and (\ref{37}) follow from (\ref{32}),
(\ref{34}), (\ref{36}) and (\ref{38}) respectively; (\ref{39})
follows from Lemma 2.7;  part~2) of the lemma follows from Lemma
2.4; part~3) of the lemma follows from an inductive application of
(\ref{2.02}); part~4) of the lemma easily follows from Lemma 2.6.
This completes the proof of the lemma.

\section{Proof of Main Result}
Let
$$
\Psi^\nu=\Phi_0\circ\Phi_1\circ\cdots\circ\Phi_\nu,~~~
\nu=0,1,\cdots
$$
Then $\Psi^\nu:D_{\nu}\times \Lambda_{\nu}\rightarrow D_0,$ and
$$
H\circ\Psi^\nu=H_{\nu}=N_{\nu}+P_{\nu},~~~ \nu=0,1,\cdots
$$
where $\Psi^0=id.$

Denote
$$
\Lambda_*=\bigcap_{\nu=0}^\infty {\Lambda}_\nu,~~~
{G}_*=D\left(\frac {r_0}2,\frac {s_0}2\right)\times {\Lambda}_*.
$$
Then $\Lambda_*$ is a Cantor-like set consisting of non-resonant
frequencies, and moreover, a measure estimate similar to that in [9]
(also [8], [10]) yields that $|\Lambda\setminus
{\Lambda}_*|=O(\gamma_0^{\frac1{n_*-1}})$.

By Lemma~3.1~2), it is easy to see that $N_\nu$ converges uniformly
on $G_*$ to
$$
N_\infty=e_\infty+\langle
\Omega_{\infty},y\rangle+\frac12\left\langle {y\choose z}, {\mathcal
M}^\infty {y\choose z}\right\rangle+h_0(x,y,z,\lambda)
$$
with
\begin{eqnarray}
|e_\infty-e_0|_{\Lambda_*}&=&O(\gamma_0^{n+1}s_0\varepsilon_0\zeta_0),\nonumber\\
|\Omega_\infty-\Omega_0|_{\Lambda_*}&=&O(\gamma_0^{n+1}s_0\varepsilon_0\zeta_0),
\nonumber\\
|\omega_\infty-\omega_0|_{\Lambda_*}&=&O(\gamma_0^{n+1}s_0\varepsilon_0\zeta_0),
\nonumber\\
|{\mathcal M}^\infty-{\mathcal M}^0|_{{\mathcal
D}(\frac{r_0}2)\times
\Lambda_*}&=&O(\gamma_0^{n+1}\varepsilon_0\zeta_0).\nonumber
\end{eqnarray}

And as in [7], we get the convergence of $\Psi^\nu$ on $G_*$ with
the estimate
$$
|\Psi^\infty-id|_{G_*}=O(\varepsilon_0\zeta_0)=O(\varepsilon_0^{\frac89}).
$$
Thus, we obtain that the perturbed system (\ref{1.1}) possesses an
analytic, quasi-periodic, invariant torus with the Diophantine toral
frequency $\omega_\infty(\lambda)=-E^\top(\lambda)\cdot
\Omega_\infty(\lambda)$  for each $\lambda\in \Lambda_*$. By
Lemma~3.1~3), we have
$$
(\Omega_{\infty}(\lambda))_{i}=(\Omega_{0}(\lambda))_i,~~~\lambda\in
{\Lambda}_*,~i=1,2,\cdots, n_0,
$$
i.e., the perturbed pseudo-frequencies  preserve the  first $n_0$
components of their corresponding ones.

In particular, when $n_0=l$, it is clear that $U^0=[A^0]$, $U=[A]$,
${\rm diag}(I_{n_0}, O)=I_l$, and, ${\rm diag}(O, I_{l-n_0})=O$.
Hence, $\Omega_{\nu}\equiv \Omega_0$ for all $\nu=0,1,\cdots$, i.e.,
$\Omega_{\infty}\equiv \Omega_0,\omega_{\infty}\equiv \omega_0$. So
we obtain that when $[A]$ is nonsingular, the Diophantine
frequencies remain unchanged under small perturbations.

\section{Some Examples}
In this section we give some examples to illustrate our results. At
first, we give an example for the persistence of invariant tori in
generalized Hamiltonian systems.

{\bf Example 1.} We consider the following unperturbed system:
$$N(y,u)=y+\frac{1}{2}y^{2}+\frac{1}{2}(u^{2}-v^{2}),$$
where $y,u,v\in R^{1},x=(x_1,x_2)^\top\in T^2,$ that is,
$l=2,n=1,m=1,$ i.e., the system is an odd dimensional generalized
Hamiltonian. The structure matrix in tangent direction $I$ is
assumed to be
$$
I={ \left(
\begin{array}{ccc}
   0 & \alpha & \beta\\
   -\alpha & 0 & -\gamma\\
   -\beta & \gamma & 0
\end{array}
\right) },
$$
where $\alpha,\beta,\gamma$ are arbitrary real numbers with
$|\alpha|+|\beta| +|\gamma|\neq 0$. It is easy to see that
\begin{align*}
\Omega=1+y,~~ \omega={ \left(
\begin{array}{c}
   -\alpha(1+y)\\
   -\beta(1+y)
\end{array}
\right) },~~ A=(1).
\end{align*}
It is easy to verify that the R\"{u}ssmann condition is not
satisfied, but $A$ is always nonsingular. So by Theorem 1.1 we
obtain that the majority 2-tori will persist with unchanged toral
frequency.

Then we give an example to illustrate the persistence of invariant
tori on sub-manifolds in generalized Hamiltonian systems. For the
persistence of elliptic invariant tori and mixed type of invariant
tori on sub-manifolds, the reader can refer to [11] for details.

{\bf Example 2.} We consider the following unperturbed system:
$$N(y,u)=\frac{1}{2}y_1^{2}+\frac12y_2^2+\frac{1}{2}(u_1^{2}-v_1^{2})
+\frac{1}{2}(u_2^{2}-2v_2^{2}),$$ where $u_i,v_i\in R^{1},i=1,2,$
and  $x=(x_1,x_2)^\top\in T^2,y=(y_1,y_2,y_3)^\top\in R^3$, that is
$l=3,n=2,m=2,$ i.e., the system is an odd dimensional generalized
Hamiltonian. The structure matrix in tangent direction $I$ is
assumed to be $$ I={ \left(
\begin{array}{ccccc}
   0 & 0 & 0 & -1& 0\\
   0 & 0 & 0 & 0& -1\\
   0 & 0 & 0 & 0& 0\\
  1 & 0 & 0 & 0& 1\\
  0 & 1 & 0 &-1& 0
\end{array}
\right) }.
$$
We consider the persistence of invariant tori on sub-manifold
$M:y_3=a,a\in R$. It is easy to see that
\begin{align*}
\Omega={ \left(
\begin{array}{c}
   y_1\\
   y_2\\
   0
\end{array}
\right) },~~ \omega={ \left(
\begin{array}{c}
   y_1\\
   y_2
\end{array}
\right) },~~ A={ \left(
\begin{array}{ccc}
  1 & 0 & 0\\
  0 & 1 & 0\\
  0 & 0 & 0
\end{array}
\right) }.
\end{align*}
By simple verification we see that the R\"{u}ssmann condition holds
and $A$ is always singular on the sub-manifold $M$. So by Theorem
1.1 we have that the first two components of $\Omega$ remain
unchanged. And by the form of $\omega(\lambda)$ we obtain that the
majority 2-tori will persist with unchanged toral frequency.
\\[3mm]
{\bf Acknowledgements} The authors express their sincere thanks to
Professor Yong Li for his instructions and encouragement.

\end{document}